\documentclass[preprint,1p,times,english]{elsarticle}
\usepackage[all]{xy} 
\usepackage{latexsym,amsmath,amssymb}
\usepackage{babel}
\usepackage{verbatim}
\usepackage{graphicx}
\usepackage{caption}
\usepackage{subcaption}

\begin{document}
\newtheorem{teorema}{Theorem}[section]
\newtheorem{lemma}[teorema]{Lemma}
\newtheorem{proposizione}[teorema]{Proposition}
\newtheorem{corollario}[teorema]{Corollary}
\newtheorem{osservazione}[teorema]{Remark}
\newtheorem{remark}[teorema]{Remark}
\newtheorem{remarks}[teorema]{Remarks}
\newtheorem{definizione}[teorema]{Definition}
\newtheorem{condizione}[teorema]{Condition}
\newtheorem{notazioni}[teorema]{Notations}
\newtheorem{examples}[teorema]{Examples}
\newtheorem{example}[teorema]{Example}
\newtheorem{definizioni}[teorema]{Definitions}
\newtheorem{definitions}[teorema]{Definitions}
\newtheorem{assunzioni}[teorema]{Assumptions}
\newcommand{\erre}{\mbox{$\mathbb{R}$}}
\newcommand{\enne}{\mbox{$\mathbb{N}$}}
\newcommand{\bvi}{\bigvee_{i=1}^{\infty} a_{i, \varphi(i)}}
\newcommand{\bvib}{\bigvee_{i=1}^{\infty} b_{i, \varphi(i)}}
\newcommand{\bvic}{\bigvee_{i=1}^{\infty} c_{i, \varphi(i)}}
\newcommand{\bvid}{\bigvee_{i=1}^{\infty} d_{i, \varphi(i)}}
\newcommand{\bvie}{\bigvee_{i=1}^{\infty} e_{i, \varphi(i)}}
\newcommand{\bviA}{\bigvee_{i=1}^{\infty} A_{i, \varphi(i)}}
\newcommand{\bviB}{\bigvee_{i=1}^{\infty} B_{i, \varphi(i)}}
\newcommand{\Fs}{\mbox{$\cal F$}}
\newcommand{\Ps}{\mbox{$\cal P$}}
\newcommand{\Gs}{\mbox{$\cal G$}}
\newcommand{\As}{\mbox{$\cal A$}}
\newcommand{\Bs}{\mbox{$\cal B$}}
\newcommand{\ttilde}{~}
\newcommand{\Ffreccia}{\mbox{$ \stackrel{\,\,\,\mathcal F\,\,\,}{\rightarrow}$}}
\newcommand{\uFfreccia}{\mbox{$ \stackrel{\,\,\,u_{\mathcal F}\,\,\,}{\rightarrow}$}}
\newcommand{\rFfreccia}{\mbox{$ \stackrel{\,\,\,r_{\mathcal F}\,\,\,}{\rightarrow}$}}
\newcommand{\oFfreccia}{\mbox{$ \stackrel{\,\,\,o_{\mathcal F}\,\,\,}{\rightarrow}$}}
\newcommand{\DFfreccia}{\mbox{$ \stackrel{\,\,\,D_{\mathcal F}\,\,\,}{\rightarrow}$}}
\newcommand{\sFfreccia}{\mbox{$ \stackrel{\,\,\, \sigma_{\mathcal F}\,\,\,}{\rightarrow}$}}

\begin{frontmatter}
\title{\bf Vitali-type theorems for filter convergence related to  vector lattice-valued modulars  
and applications to stochastic processes\tnoteref{supporto}}
\tnotetext[supporto]{ Supported by University of Perugia and G. N. A. M. P. A.}

\author[]{A. Boccuto} \ead{boccuto@dmi.unipg.it}
\author[]{D. Candeloro \corref{cor1}} \ead{candelor@dmi.unipg.it}
\cortext[cor1]{corresponding author address: Department of Mathematics and Computer Sciences 1, Via Vanvitelli 06124 Perugia (Italy)}

\author[]{A. R. Sambucini} \ead{matears1@unipg.it}

\address{Department of Mathematics and Computer Sciences - 06124 Perugia (Italy)}
\date{May 23th, 2013}

\begin{abstract} 
A Vitali-type theorem for vector lattice-valued modulars
with respect to filter con\-ver\-gen\-ce is proved. Some applications are given
 to modular con\-ver\-gen\-ce theorems for moment operators
in the vector lattice   setting, and also for the Brownian motion and related stochastic processes.
\end{abstract}

\begin{keyword}
 Vector lattice \sep   filter convergence \sep   modular    \sep     Vitali convergence theorem \sep   moment operator \sep 
   Brownian motion \sep It\^{o} integral\\ 
\MSC[2010] 28B15 \sep 41A35 \sep 46G10
\end{keyword}
\end{frontmatter}

\section*{Introduction} The Vitali theorem has been widely studied in the literature.
For a historical survey, see for instance  \cite{CHOKSI} and the bibliography therein. 
This theorem
is a well-known milestone in integration theory, directly related with $L^p$ convergence
and weak compactness. However it has also many applications in several branches of Mathematics, for example in
reconstruction of signals, integral and discrete operators (see \cite{bm1, bm2006, BMV, BBV1, BBV2, BBV3, VZ, CV}).
The Vitali theorem has been investigated even in  the context of vector lattices with
respect to different types of convergence (see  \cite{BC} and its bibliography) and in abstract 
convergence groups, not necessarily endowed with an order structure 
(see \cite{bccsmb}).

Also the modular spaces (see for instance  \cite{BMV, KOZLOWSKI}) are a rich field of  research.
They are a natural generalization of the Lebesgue spaces $L^p$ and contain
as particular cases the Orlicz and the Musielak-Orlicz spaces (see e.g. \cite{MUSIELAK, RAO}).

In modular spaces the Vitali theorem is often used in
the problem of approximating a real-valued  function $f$ by means of Urysohn-type 
or sampling operators $(T_n f)_n$, (see also \cite{bm1, bm2006, BV1, BV2, BV3, 
VINTI, AV, AV1,AV2, bd1, bd2, bbdm,VZ1}).
These operators are particularly useful in order to  approximate a continuous or analog signal by 
means of discrete samples, and therefore they are  widely applied for example in 
reconstructing images and videos.

In the classical literature, a signal $f$ is  considered as a real-valued function, defined on a
(bounded or not) time interval. In these situations, 
usually it happens that the original signal or image is unknown, and 
it must be reconstructed from some sampled
values of the type $f\Bigl(\frac{k}{w}\Bigr)$, where $k \in \mathbb{Z}$, and $w$ is  a fixed rate.
So it would be advisable to consider the studied signal, for example, as a map defined on a (finite or not)
time interval and taking values in a suitable space  of random variables (for more details related
with these topics see also \cite{BMV} and the literature therein). 
In this context, it is natural to consider  the tool of \textit{vector lattices}.
Moreover, also domains with infinite measure, like for instance an infinite time interval, can be  treated.

In this paper we consider vector lattices with convergences
generated by filters (see also \cite{bdpschurfilters}), 
which generalize classical ones.
For example, when the filter involved is the class of all subsets  of $\mathbb{N}$
having asymptotic density one, we obtain the so-called statistical convergence.
We introduce the theory of modular spaces in the vector lattice  context:
integration theory can be considered as a 	particular case of the theory of modulars, even in 
our setting, so we give a general version of the Vitali and Lebesgue dominated convergence 
theorems for modular convergence, with  applications to  integrals
for vector lattice-valued functions, defined  on sets of possibly infinite measure.
In Section 4  we deal with 
moment kernels in the vector lattice  setting,
in order to construct a suitable sequence of convolution type operators; moreover
in Subsection \ref{bm}, among the applications, 
we give convergence results for the It\^{o} integral with respect  to Brownian motion.
Finally, the last Section is devoted to the proof of the main theoretical result (Theorem \ref{unico}) of Section 4.
\section{Preliminaries}
We begin with recalling some basic properties of vector lattices and filter convergence.
\subsection{Vector lattices}
A vector lattice (also called Riesz space)  ${\bf X}$ is said to be \textit{Dedekind complete} iff
every nonempty subset
$B$ of ${\bf X}$, bounded from above, has a lattice supremum  in ${\bf X}$, denoted by $\bigvee B$.
\\
From now on ${\bf X}$ denotes a Dedekind complete vector lattice,  ${\bf X}^+$ is the set of all strictly positive elements of
${\bf X}$, and we set ${\bf X}^+_0={\bf X}^+ \cup \{0\}$. For all $x \in {\bf X}$, let $|x| := x \bigvee (-x)$.
We add to ${\bf X}$ an extra element $+\infty$, extending order and operations in a natural way, set
 $\overline{{\bf X}}={\bf X} \cup \{+ \infty\}$, $\overline{{\bf X}^+_0}={\bf X}^+_0 \cup \{+ \infty\}$ and assume by convention
$0 \cdot (+ \infty) = 0$.
\\
A sequence $(p_n)_n$ in ${\bf X}$ is called \textit{$(o)$-sequence} iff it is decreasing and $\wedge_n p_n=0$.
\\
A \textit{strong unit} of ${\bf X}$ is an element $e$, such that for every $x \in {\bf X}$ there is a positive real number
$c$ with $|x|\leq c \, e$.
\\
Let $ \bf{X}_1$, $ \bf{X}_2$, ${\bf X}$ be three Dedekind complete vector lattices.
We say that $( \bf{X}_1, \bf{X}_2, {\bf X})$ is a \em product triple \rm iff a map $\cdot : \bf{X}_1
\times  \bf{X}_2 \to {\bf X}$ (\em product\rm) is given, satisfying  natural conditions of compatibility,
see \cite[Assumption 2.1]{BC2011} and \cite{bms}. From now on we suppose that:
\begin{itemize}
\item[ \textbf{($H_0$)}]
 $( \bf{X}_1,  \bf{X}_2, {\bf X})$ is a product triple,  where ${\bf X},  \bf{X}_1$ are endowed 
with strong units $e$, $e_1$ re\-spe\-ctively.
\end{itemize}
\subsection{Filter convergence}
Let $Z$ be any fixed countable set. A class  
$\mathcal{F}$ of subsets of $Z$ is called  a \textit{filter} of $Z$ iff $\emptyset 
\not \in {\mathcal F}$,  $A \cap B \in {\mathcal F}$ whenever $A,B \in 
\mathcal{F}$ and for each $A \in \mathcal{F}$ and $B \supset A$ we get  $B\in \mathcal{F}$.
We denote by ${\mathcal F}_{\text{cofin}}$ the  filter of all cofinite subsets of $Z$.
A filter of $Z$ is said to be \textit{free} iff it  contains ${\mathcal F}_{\text{cofin}}$.
An example of free filter (different from ${\mathcal F}_{\text{cofin}}$) is the family ${\mathcal F}_{\text{d}}$ of all subsets of $\mathbb{N}$,
having asymptotic density one. 
In this paper we only deal with free filters.\\

We now give the notion of filter convergence in the  vector lattice  setting (see also 
\cite{bdpschurfilters}).
\begin{definizione}\label{filterconv} \rm
Let ${\mathcal F}$ be any filter of $Z$.
A sequence $(x_z)_{z \in Z}$ in ${\bf X}$ 
\textit{$(o_{\mathcal F})$-converges to $x \in {\bf X}$}    ($x_z \oFfreccia x$) iff there exists 
an $(o)$-sequence $(\sigma_p)_p $ in ${\bf X}$  such
that for all $p \in \mathbb{N}$ the set $\{ z \in Z:  |x_z - x| \leq \sigma_p  \}$
is an element of ${\mathcal F}$.\\
A sequence $(x_z)_{z \in Z}$ in ${\bf X}$ \textit{$(r_{\mathcal F})$-converges to $x
\in {\bf X}$}    ($x_z \rFfreccia x$) iff there exist a positive element $u \in {\bf X}$ and  an $(o)$-sequence $(\varepsilon_p)_p \in 
\mathbb{R}^+$  such
that for all $p \in \mathbb{N}$ the set $\{ z \in Z: |x_z - x| \leq  \varepsilon_p \, u\}$ is an element of ${\mathcal F}$.\\

Clearly a sequence $(x_z)_z$ in ${\bf X}$ \textit{$(o)$}-converges to $x \in {\bf X}$ (in the classical sense) iff
it $(o_{{\mathcal F}_{\text{cofin}}})$-converges to $x$, and similarly with $(r)$-convergence.
\end{definizione}
\section{Convergence theorems}

Let $G$ be any infinite set, $\Sigma\subset {\mathcal P}(G)$ be a 
algebra and $\mu:\Sigma \to \overline{( {\bf X_2})^+_0}$ be a finitely additive measure.
We shall denote by 

$\Sigma_b$ the family of all sets $B\in \Sigma$ with $\mu(B)\in  \bf{X}_2$.
For every $A \in \Sigma$ and $r \in  \bf{X}_1$, 
we shall denote by $r\cdot 1_A : G \rightarrow  \bf{X}_1$ the function whose values $r\cdot 1_A(t)$ are $r$ when $t\in A$ and 0 otherwise.
\subsection{Modulars}
We now introduce the modulars in the setting of vector lattices
(for the classical case and a related literature, see e.g. \cite{BMV, KOZLOWSKI}).
Let $T$ be a linear sublattice of $ \bf{X}_1^G$, such that 
\begin{itemize}
\item $e_1 \cdot 1_A \in T$ for every $A \in \Sigma_b$;
\item if $f \in T$ and $A \in \Sigma$, then $f \cdot 1_A \in T$.
\end{itemize} 
\noindent
A functional $\rho:T \to \overline{{\bf X}^+_0}$ is said to be a
\textit{modular} on $T$ iff it satisfies the following properties:
\begin{description}
\item[($\rho_0$)] $\rho(0)=0$;
\item[($\rho_1$)] $\rho(-f)=\rho(f)$ for every $f \in T$;
\item[($\rho_{2}$)] $\rho(\alpha_1 f + \alpha_2 h) \leq \rho(f) + \rho(h)$ for every
$f$, $h \in T$ and for any non negative real numbers $\alpha_1, \alpha_2$\,  with
$\alpha_1 + \alpha_2 =1$.
\end{description}
We say that:
\begin{description}
\item[($\rho_m$)  ] A modular $\rho$ is  \textit{monotone} iff $\rho(f) \leq \rho(h)$ for every $f$, $h \in T$ with
$|f|\leq |h|$. In this case, if $f \in T$, then $|f| \in T$ and hence $\rho(f) = \rho(|f|)$.
\item[($\rho_f$)  ]
A modular $\rho$ is \textit{finite} iff  for every $A\in \Sigma_b$  and every 
$(o)$-sequence $(\varepsilon_p)_p$ in $\mathbb{R}^+$, the sequence $(\rho(e_1 \ 
\varepsilon_p\ 1_A))_p$ is $(r)$-convergent to $0$.
(For the case ${\bf X} = \bf{X}_1=  \bf{X}_2=\mathbb{R}$, see also \cite{BMV}).
\end{description}
We now prove the following
\begin{proposizione}\label{propertyglobalfiniteness}
Let $( \bf{X}_1, \bf{X}_2, \bf{X})$ 
be a product triple  satisfying \textbf{($H_0$)}, and assume that 
$\rho$ is a  finite modular, $(B_m)_m$ is a sequence
in $\Sigma_b$ and $(\varepsilon_p)_p$ is any  $(o)$-sequence in $\mathbb{R}$. Then there 
exists a strictly increasing  mapping $m\mapsto p(m)$ from $\mathbb{N}$ to 
$\mathbb{N}$, with
$$\rho(e_1 \ \varepsilon_{p(m)} 1_{B_m})\leq  \frac{e}{m} \quad \text{for every  }m \in \mathbb{N}.$$
\end{proposizione}
{\bf Proof:}
Fix a sequence $(B_m)_m$  in $\Sigma_b$, and any $(o)$-sequence 
$(\varepsilon_p)_p$  in $\erre$. By finiteness of $\rho$, we get
$$(r)\lim_{p \to +\infty} \rho(e_1 \, \varepsilon_p \, 1_{B_m})=0 \quad
\text{for every  } m \in \enne.$$
This means that  a sequence $(u_m)_m$ exists in   ${\bf X}^+$ with the following property: for every
$m$, $p\in \mathbb{N}$ there is $\lambda(m,p) \in \mathbb{R}^+$ such 
that $(\lambda(m,p))_p$ is an $(o)$-sequence for all $m$, and
\begin{eqnarray}\label{unoo}
0 \leq  \rho(e_1 \  \varepsilon_p \, 1_{B_m}) \leq  \lambda(m,p)u_m\,\quad \text{for all  }m, p.
\end{eqnarray}
For every $m \in \mathbb{N}$ there is a positive real number $h_m$ such that $u_m\leq h_m\ e$.
Therefore, from (\ref{unoo}) we also have
\begin{eqnarray}\label{duee}
0 \leq \rho(e_1 \ \varepsilon_p \, 1_{B_m}) \leq   h_m\lambda(m,p)\, e \quad
\text{whenever}\,\,  m,p\in \enne.
\end{eqnarray} 
Now, for each integer $m$, a corresponding integer 
$p(m)$ can be found, such that the sequence $(p(m))_m$ is strictly increasing and  
\begin{eqnarray}\label{tree}
h_m\lambda(m,p(m))\leq \frac{1}{m}  \quad \text{for all  } m.
\end{eqnarray} 
The assertion follows from (\ref{duee}) and (\ref{tree}). $\quad \Box$
\\

We now give the concept of (equi)-absolute  continuity in the context of modulars and filter 
convergence: this is one of the crucial tools in the Vitali theorem (for a classical modular 
version, see \cite[Theorem 2.1]{BMV}). 
\begin{description}
\item[($a_{\rho}$)]
We say that $f \in T$ is \textit{$(o_{\mathcal F})$-absolutely
continuous with respect to the modular $\rho$} (or in short \textit{absolutely continuous}) iff there is a positive real constant $\alpha$,
satisfying the following properties:
\item[($a_{\rho}(1)$)] for each $(o)$-sequence $(\sigma_p)_p$ in $  \bf{X}_2^+$ there exists an $(o)$-sequence $(w_p)_p$ in ${\bf X}^+$
such that for all $p \in \enne$ and whenever $\mu(B) \leq \sigma_p$ we get
$\rho(\alpha f \, 1_B) \leq w_p$;
\item[($a_{\rho}(2)$)] there is an $(o)$-sequence $(z_m)_m$ in ${\bf X}^+$ such that to
each $m \in \enne$ there corresponds a set $B_m \in \Sigma_b$
with $\rho(\alpha f \, 1_{G \setminus B_m}) \leq z_m$.
\item[($ac_{\rho}$)]
Given a modular $\rho$ and any free filter ${\mathcal F}$ of $Z$,
we say that a sequence $f_z:G \to {\bf X_1}$, $z \in Z$, is
\textit{$\rho$-${\mathcal F}$-equi-absolutely continuous}, or in short
\textit{equi-absolutely continuous,} iff
there is  $\alpha \in \mathbb{R}^+$, satisfying the following two conditions:
\begin{description}
\item[($ac_{\rho}(1)$)] for every $(o)$-sequence $(\sigma_p)_p$ in 
$ \bf{X}_2^+$  there are an $(o)$-sequence $(w_p)_p$ in ${\bf X}^+$
and a sequence $(\Xi^p)_p$ in ${\mathcal F}$,
with $\rho(\alpha f_z \, 1_B) \leq w_p$
whenever $z \in \Xi^p$ and $\mu(B) \leq 
\sigma_p$, $p\in \enne$;
\item[($ac_{\rho}(2)$)] there are an  $(o)$-sequence $(r_m)_m$ in ${\bf X}^+$ 
and a sequence $(B_m)_m \in \Sigma_b$ such that
\begin{eqnarray}\label{equiac2}
\Lambda^{m}:=\{z \in Z: \rho(\alpha f_z \, 1_{G  \setminus B_m}) \leq r_m \} \in {\mathcal F}
\quad \text{for all  } m\in \mathbb{N}.
\end{eqnarray} 
\end{description}
\end{description}
In \cite{bm1} some sufficient conditions are given  for (equi)-absolute continuity with respect to modulars.
\subsection{A Vitali-type theorem}
As we already pointed out in the Introduction, the Vitali theorem has a fundamental importance in 
Approximation Theory, also when $\mu(G)=+\infty$: several applications  can be found in \cite{bm2006, 
BMV, BUTZER, HOLMES, MUSIELAK, VINTI} and the literature therein. Together with equi-absolute 
continuity, also some notions of convergence are  needed, so we recall the concepts of filter uniform 
convergence and convergence in measure (see also  \cite{BC}).
Let ${\mathcal F}$ be any fixed free filter of $Z$.
\begin{definizioni} \rm \label{u-misura}
\ttilde
\begin{description} \vskip-1cm
\item[(\ref{u-misura}.1)]  A sequence of  functions 
$(f_z)_{z \in Z}$ in $ \bf{X}_1^G$ {\em  \textit{ $(r_{\mathcal F})$}-converges uniformly} (shortly 
\textit{converges uniformly}) to $f$,
iff there exists an $(o)$-sequence $(\varepsilon_p)_p$ in $\mathbb{R}^+$ with
\[ \{ z \in Z : \bigvee_{t \in G} \, |f_z(t)-f(t)| < \varepsilon_p \, e_1 \} \in {\mathcal F}
\quad \text{for  any  } p \in \mathbb{N}.\] In this case we write:
$$(r_{\mathcal F})\lim_z \Bigl( \bigvee_{t \in G} \, |f_z(t)-f(t)|\Bigr)=0.$$
\item[({\bf\ref{u-misura}.2)}] Given a sequence  $(f_z)_z$ in $ \bf{X}_1^G$ and $f \in  \bf{X}_1^G$, 
we say that
$(f_z)_z$ \textit{$(r_{\mathcal F})$-converges in measure}
(shortly \textit{$\mu$-converges}) to $f$,
iff there are  two $(o)$-sequences,
$(\varepsilon_p)_p$ in $\erre^+$, $(\sigma_p)_p$  in $ \bf{X}_2^+$,
and a double sequence $(A_z^p)_{(z,p) \in Z 
\times \mathbb{N}}$ in $\Sigma$ such that
\begin{itemize}
\item  $A_z^p \supset \{t \in G: |f_z(t)-f(t)| \not\leq \varepsilon_p \ e_1\}$ for every $z \in Z$ 
and $p \in \enne$, and 
\item $\{z \in Z:  \mu(A_z^p) \leq \sigma_p \} \in {\mathcal F}$  for all $p \in \mathbb{N}$. 
\end{itemize}
\end{description}
\end{definizioni}
It is easy to see that {(\bf\ref{u-misura}.1)}  implies {(\bf\ref{u-misura}.2)}. 
If $(f_z)_z, (h_z)_z$ $\mu$-converge to $f,h$ 
respectively, then $(f_z \pm h_z)_z$ $\mu$-converges to $f \pm h$ and
$(|f_z|)_z$ $\mu$-converges to $|f|$, (see also \cite{BC}).
Moreover,  when $ \bf{X}_1= \bf{X}_2=\mathbb{R}$ and the convergence involved is the usual 
one, these notions coincide with the classical ones.\\

We now prove a Vitali-type theorem, which links the theories of modulars and vector lattices, in the context of $(o_{\mathcal F})$-convergence.
\begin{teorema}\label{trediecibis} \mbox{\rm (\textbf{Vitali})}
Let $\rho$ be a monotone and  finite modular,   and 
$\mathcal F$ be a fixed free filter of $Z$. If $(f_z)_z$ is a sequence in $T$, 
$\mu$-convergent  to $0$ and  equi-absolutely continuous, then there is a positive number $\alpha$ with
$$(o_{\mathcal F})\lim_z \rho(\alpha \, f_z)=0.$$
\end{teorema}
{\bf Proof:}
Let $\alpha \in \mathbb{R}^+$ be related with equi-absolute continuity, and
$(r_m)_m$ in ${\bf X}^+$ be an $(o)$-sequence associated with ($ac_{\rho}(2)$); let moreover
$(\sigma_p)_p$ in $ \bf{X}_2^+$,  $(\varepsilon_p)_p$ in $\mathbb{R}^+$ be 
related with $\mu$-convergence of the sequence $(f_z)_z$ to $0$, and
$(w_p)_p$ in ${\bf X}^+$ be associated with $(\sigma_p)_p$ and ($ac_{\rho}(1)$).

Choose any sequence $(B_m)_m$ in $\Sigma_b$, and any corresponding sequence 
$(\Lambda^{m})_m$ in $\mathcal{F}$, satisfying (\ref{equiac2}).
Without loss of generality, we assume that   the $f_z$'s are positive. 
For every $z\in Z$, $p \in \mathbb{N}$, let 
$A_z^p \in \Sigma$ be such that 
$A_z^p \supset \{t \in G: f_z(t) \not\leq \varepsilon_p \ e_1  \}$,  according to  $\mu$-convergence.
For every $z\in Z$,  $m$, $p \in \mathbb{N}$ and  $t \in B_m$ we have
\begin{align}\label{medicinaa}
\frac{\alpha}{3} f_z(t)=\frac1{3} \Bigl(\alpha f_z(t) 1_{G \setminus B_m}(t) + \alpha f_z(t)
1_{B_m \setminus A_z^p}(t) + \alpha f_z(t)  1_{B_m \cap A_z^p}(t) \Bigr).
\end{align}
By the properties of the modular $\rho$  we obtain (with obvious meaning of the notations):
\begin{align}\label{matematicaa}
\rho \nonumber \Bigl(\frac{\alpha}{3} f_z \Bigr)&\leq \rho(\alpha \, f_z 1_{G \setminus B_m})+
\rho(\alpha \, f_z 1_{B_m \setminus A_z^p})+
\rho(\alpha \, f_z 1_{B_m \cap A_z^p}) = \\ &=I_1+I_2+I_3 \quad
\text{ for all } z\in Z,\ m,p \in \enne.
\end{align}
By (\ref{equiac2}) we get 
\begin{align}\label{iuno}
I_1=\rho(\alpha f_z 1_{G \setminus B_m}) \leq r_m \quad \text{ for all } z \in \Lambda^{m}.
\end{align}
Moreover, since $\mu(A_z^p)\leq \sigma_p$ by 
convergence in measure, then, thanks to property
$(ac_{\rho}(1))$ of equi-absolute continuity, 
for each $p \in \enne$ there exists an element $\Xi^{p} \in {\mathcal F}$, with
\begin{align}\label{itre}
I_3 = \rho(\alpha f_z 1_{B_m \cap A_z^p}) \leq \rho(\alpha f_z 1_{A_z^p}) \leq w_p
\quad \text{ whenever } m \in \enne \text { and } z \in \Xi^{p}.
\end{align}
In order to estimate $I_2$, 
observe that  by Proposition 
\ref{propertyglobalfiniteness} there is
a strictly increasing mapping $m\mapsto p(m)$ from $\mathbb{N}$ to $\mathbb{N}$, with
$$\rho(\alpha  \varepsilon_{p(m)} \ e_1  1_{B_m})\leq \frac{e}{m} \quad \text{for all }m;$$
now,  for all  $z\in Z$, $m \in \enne$ and $t \in  B_m$ we see that
$$0 \leq \alpha \, f_z(t) 1_{B_m \setminus  A_z^{p(m)}}(t) \leq \alpha\varepsilon_{p(m)} e_1 \  1_{B_m}(t) $$ 
and, consequently,
\begin{align}\label{idue}
I_2=\rho(\alpha f_z 1_{B_m \setminus A_z^{p(m)}})\leq \frac{e}{m}.
\end{align}
Thus from (\ref{matematicaa}), (\ref{iuno}), (\ref{itre}) and (\ref{idue}) it follows that for every $m$
there are $\Lambda^{m}$, $\Xi^{p(m)} \in {\mathcal F}$ such that
$$\rho\Bigl(\frac{\alpha}{3} f_z \Bigr) \leq \rho(\alpha \, f_z 1_{G \setminus B_m})+
\rho(\alpha \, f_z 1_{B_m \setminus A_z^{p(m)}})+
\rho(\alpha \, f_z 1_{B_m \cap A_z^{p(m)}}) \leq r_m + \frac{e}{m}  + w_{p(m)}$$
whenever $z \in \Lambda^{m} \cap\Xi^{p(m)}$,
and so, since $m\mapsto r_m + \displaystyle{\frac{e}{m}}  + w_{p(m)}$ is an $(o)$-sequence, we conclude
$$(o_{\mathcal F})\lim_z \rho\Bigl(\frac{\alpha}{3} \, f_z\Bigr)=0,$$
that is the assertion. $\quad \Box$\\

An easy consequence of Theorem \ref{trediecibis} is the following
\begin{teorema}\label{lebesgue}
\mbox{\rm (\bf Lebesgue dominated convergence theorem\rm)}
Let ${\bf X}$, $\rho$, ${\mathcal F}$ be as in Theorem  \ref{trediecibis} and assume that $(f_z)_z$ 
is a sequence in $T$,  $\mu$-convergent to $0$. Suppose that there are  an absolutely continuous
function $g$ in $T$ and an element $F_0 \in {\mathcal F}$, such that $|f_z(t)| \leq g(t)$ for all $z \in F_0$ and $t \in G$.
Then there is a positive real number $\alpha$ with $(o_{\mathcal F})\lim_z \rho(\alpha \, f_z)=0$.
\end{teorema}

\section{Examples and applications}
We give here some applications of the results 
obtained above in different fields of Mathematics, 
as integration theory, moment  operators, stochastic processes and Brownian  motion.
\subsection{Integration theory}
Let $G$, $\Sigma$ and $\mu$ be as above.
We now extend the integration theory investigated  in \cite{BC} to the case in which $\mu$ can 
assume the value $+\infty$, in the setting of $(o_{\mathcal F})$-convergence.\\
\begin{definizioni}\rm 
A function $f: G \rightarrow {\bf X_1}$ is said to  be \textit{simple} iff $f(G)$ is a finite set and
$f^{-1} (x) \in \Sigma$ for every $x \in {\bf X}_1$. The space of all simple functions 
is denoted by ${\mathcal S}$.
\end{definizioni}

Let $L^*$ be the set of all simple functions $f\in {\mathcal S}$ vanishing outside a set of $\mu$-finite measure.
If $f \in L^*$, we denote its elementary  integral with $\displaystyle{\int_G f(t) \,d\mu(t)}$, as usual.

\begin{remark}\label{additivita} \rm
It is clear that the functional $\iota: L^* \rightarrow {\bf X}$ defined as
\begin{eqnarray}\label{modulareintegrale}
\iota(f) :=\int_G |f(t)| \, d\mu(t), \quad f \in L^*,
\end{eqnarray}
is a monotone modular, and it is also linear and  additive on positive functions and constants.

Also finiteness of $\iota$ is easy to see: indeed, let $A \in \Sigma_b$ 
and $(\varepsilon_p)_p$ be any $(o)$-sequence in $\erre^+$. According to our definitions, 
$\iota(\varepsilon_p \ e_1 1_A)=\varepsilon_p  \mu(A) \ e_1$ for all $p\in \enne$: 
since $\mu(A)  e_1$ is a positive element in ${\bf X}$, it is 
obvious that the condition of finiteness is  satisfied. 
\end{remark}
\begin{definizione}\label{integrabilitainfinita} \rm 
A positive function $f \in  \bf{X}_1^G$ is
\textit{integrable} iff there exist an equi-absolutely continuous sequence
of functions $(f_{n})_{n}$ in $L^*$, $\mu$-convergent to $f$, and a map $l:\Sigma \to {\bf  X}$, with
\begin{eqnarray}\label{zvezda} (o_{\mathcal  F})\lim_n \bigvee_{A \in \Sigma}
\Bigl|\int_A f_n(t)\, d\mu(t)-l(A)\Bigr|=0
\end{eqnarray}
(the sequence $(f_n)_n$ is said to be  \textit{defining}). Note that, in this case, we get
 $$l(A)=(o_{\mathcal F})\lim_n \int_A \,  f_n(t)\,d\mu(t) \quad \text{uniformly with respect to  } A \in  \Sigma.$$

If $f \in  \bf{X}_1^G$,  then we say that $f$ is \textit{integrable} iff $f^+$ and
$f^-$ are integrable, where $f^+(t)=f(t) \vee 0$, $f^-(t)=(-f(t)) \vee 0$, $t \in G$.
In this case, we set
\begin{eqnarray}\label{pm2}
\int_A f(t) \, d\mu(t):=\int_A f^+(t) \, d\mu(t) - \int_A f^-(t) \, d\mu(t), \quad A \in \Sigma.
\end{eqnarray}
\end{definizione}
Similarly as in \cite[Proposition 3.11]{BC}, we now prove that the integral in Definition \ref{integrabilitainfinita} is well-defined.
\begin{teorema}\label{int-inf}
In the setting described above, let  $l$ be as in (\ref{zvezda}).
Then the quantity $l(A)$ is independent of the choice of the defining sequence $(f_n)_n$.
\end{teorema}
{\bf Proof:} Let $(f_n^j)_n$, $j=1,2$, be two  defining sequences for $f$, and put
\begin{eqnarray*}
l_j(A)&=&(o_{\mathcal F})\lim_n \int_{A} \,  f_n^j(t)\,d\mu(t), \quad A \in \Sigma, \, j=1,2;
\\ h_{n}(t)&=&|f_{n}^{1}(t)-f_{n}^{2}(t)|, \, t  \in G,  \, n \in \mathbb{N}.
\end{eqnarray*}
Since $(f_{n}^{j})_n$, $j=1,2$, are equi-absolutely continuous, then
$(h_{n})_n$ is too. As $(f_{n}^{1})_n$ and $(f_{n}^{2})_n$ $\mu$-converge  to $f$, 
$(h_{n})_n$ $\mu$-converges  to $0$. 
Thus, by Theorem \ref{trediecibis} applied
to the modular $\iota: L^* \rightarrow {\bf X}$ defined in Remark \ref{additivita}, we obtain:
$$\displaystyle{(o_{\mathcal F})\lim_n\int_G h_n(t) \, d\mu(t) =0}.$$
Therefore for all $A \in \Sigma$ and $n\in \enne$ we get:
\begin{eqnarray*}
|l_{1}(A)-l_{2}(A)| &\leq&
\left|\int_{A}f_{n}^{1}(t) \,
d\mu(t)-l_{1}(A)\right| + \left|l_{2}(A)- \int_{A}f_{n}^{2}(t) \,
d\mu(t)\right|+ \int_G h_n(t) \, d\mu(t).
\end{eqnarray*}
By letting $n$ tend to infinity, we obtain $|l_{1}(A)-l_{2}(A)|=
0$, namely $l_{1}(A)=l_{2}(A)$ for all $A \in \Sigma$. $\quad \Box$
\\

\begin{remark}\label{additivitaestesa} \rm
Let $L$ be the space of all
integrable functions in the sense of Definition \ref{integrabilitainfinita}.
It is clear that in this way the elementary integral, defined in $L^*$, can be uniquely extended to a linear monotone functional $\tau: L \rightarrow {\bf X}$ defined on the whole of $L$.
In particular, it turns out that, fixed a measurable  subset $B\subset G$ and a mapping $f\in  
\bf{X}_1^G$, then $f$ is integrable in the whole  of $G$ iff both mappings $f\ 1_B$ and $f\ 
1_{G\setminus B}$ are (and of course the integral  is additive in this case). 
\end{remark}

\begin{remark}\label{perdopo} \rm
If $f:  G \rightarrow  {\bf X}$ is defined by $f(t) = h(t) \ u$, where
$u \in {\bf X}^+$ is fixed and $h: G \rightarrow \mathbb{R}$ is  Lebesgue integrable, then $f$ belongs to $L(\lambda)$
and $\tau(f) = (\int_G h d\mu) u$.
\end{remark}
\subsection{The moment operator}
Let $ \bf{X}_1={\bf X}$ be a vector lattice  with a  strong unit $e$, $ \bf{X}_2 = \mathbb{R}$,
$G = \mathbb{R}^+$, ${\mathcal B}$ be the
$\sigma$-algebra of all measurable subsets of $ \mathbb{R}^+$,  $ \lambda:{\mathcal B} \to \mathbb{R}^+$
be the Lebesgue measure and  ${\mathcal F}={\mathcal F}_{\text{cofin}}$.

 From now on we consider functions $f:  
\mathbb{R}^+ \to{\bf X}$, and the so-called \textit{moment kernels},
$M_n: \mathbb{R}^+ \to \mathbb{R}$, defined by
$$M_n(w)=w^n \cdot n \cdot 1_{]0,1[} (w),\,\, n \in \mathbb{N}, w \in  \mathbb{R}^+.$$

Let  $\widehat{L}(\lambda)$ be the set defined by
\begin{eqnarray}\label{integrale-cappuccio}
\widehat{L}(\lambda) := \{ f \in {\bf X}^{ \mathbb{R}^+}: \,f \in L(\lambda)
  \text{ and }
\displaystyle{M_n \Bigl( \frac{\cdot}{s} \Bigr)  f(\cdot) \in L(\mu)}\
\forall s \in  \mathbb{R}^+, n \in {\mathbb N} \}
\end{eqnarray}
where $\displaystyle{\mu(B):=\int_B \, \frac{dt}
{t}}$ for every $B \in {\mathcal B}$.
We consider the modular $\displaystyle{\rho (\cdot) = \int_{\mathbb{R}^+} |\cdot| ds}$.\\

In order to define and study the moment operators we need the following notions and results.
\begin{definizione}\label{order-uc}\rm
 We say that $f:\mathbb{R}^+\to{\bf X}$ 
is \textit{uniformly continuous}  iff there are a positive element $u\in {\bf X}$ and two $(o)$-sequences
$(\sigma_p)_p$ and  $(\delta_p)_p$ in $\mathbb{R}^+$ such that
$|f(x_1)-f(x_2)|\leq \sigma_p\ u$ for all $x_1$, $x_2 \in  \mathbb{R}^+$  and $p \in {\mathbb N}$ satisfying
$|x_1-x_2| < \delta_p$.
\end{definizione}

We remark here that any uniformly continuous  mapping $f$ is {\em locally bounded}, i.e. the set 
$\{|f(x)|: x\in (a,b)\}$ is bounded in ${\bf X}$, for  any bounded interval $(a,b)\subset \erre^+$. 
\begin{definizione}\label{order-lip}\rm
We say that $f:\mathbb{R}^+\to{\bf X}$ is \textit{Lipschitz}, iff there is $K \in {\bf X}^+$
with $|f(x_1)-f(x_2)|\leq K |x_1 - x_2|$ for all $x_1$, $x_2 \in  \mathbb{R}^+$.  If $f$ is 
Lipschitz, then $f$ is uniformly continuous. 
\end{definizione}
\begin{definizione} \rm
Let $I \subset \mathbb{R}$ be a connected set, and
$f:I \to {\bf X}$. We say that $f$ is \textit{uniformly differentiable on $I$} iff there exist
a bounded function $g:I \to {\bf X}$ and
two $(o)$-sequences, $(\sigma_p)_p$  and  $(\delta_p)_p$ in $\mathbb{R}^+$, with
$$\bigvee \left\{\left|\frac{f(v)-f(u)}{v-u}-g(x)\right|:
u \leq x \leq v , 0 < v-u \leq \delta_p \right\} \leq \sigma_p\ e$$
for every $p \in \mathbb{N}$. In this case we say that $g$ is the \textit{uniform derivative} of $f$ in $I$.
\end{definizione}
Note that in general the concept of uniform   derivative is a strict strengthening of the classical 
one. However if $f \in C^1(I,\mathbb{R})$ in the  classical sense, then $f$ is uniformly differentiable 
and $g= f^{\prime}$.
\\

We see that uniform differentiability and uniform  continuity in the vector lattice  setting
satisfy analogous properties as the corresponding classical ones.
\begin{proposizione} \label{uc-prodotto}
Let $f : [a,b] \to {\bf X}$, $g: [a,b] \to \mathbb{R}$. Then we have:
\begin{description}
\item[(\ref{uc-prodotto}.1)] if $f,g$ are  uniformly continuous,  then $f \cdot g$ is uniformly continuous too;
\item[(\ref{uc-prodotto}.2)] if $f$ is uniformly differentiable, then $f$ is uniformly continuous;
  \item[(\ref{uc-prodotto}.3)] if $f,g$  are  uniformly differentiable,
with uniform derivatives $f^{\prime}$, $g^{\prime}$ respectively, then $f \cdot g$ is uniformly
differentiable too, and its uniform derivative is $f^{\prime} \, g + g^{\prime} \, f$;
\item[(\ref{uc-prodotto}.4)] if $f$ is uniformly 
continuous, then $f$ is bounded and  belongs to $L(\lambda)$.
Moreover its integral function $F$ is uniformly differentiable in $[a,b]$, with uniform derivative $f$.
\end{description}
\end{proposizione}
{\bf Proof:} The proofs of statements (\ref{uc-prodotto}.$j$)   for $j=1,2,3$ are straightforward.
We prove now (\ref{uc-prodotto}.4).\\
Boundedness has already been pointed out. So we prove just integrability.
Choose $(\sigma_p)_p$ and $(\delta_p)_p$ in $\mathbb{R}^+$ according to uniform continuity of $f$.
 For each $p \in \mathbb{N}$, let $D_p:=\{ a=x_0 < x_1 < \ldots < x_{q_p}=b \}$ be a division
of $[a,b]$, with mesh $\delta(D_p) \leq \delta_p$, and take
$D_{p+1}$ in such a way that $D_{p+1}$ is a refinement of $D_p$. 
Now, set
$$f_p(t):=\sum_{i=0}^{q_p-1} \, \bigwedge_{t \in [x_i, x_{i+1}[} f(t) \cdot 1_{[x_i, x_{i+1}[}(t).$$
So, there exists a suitable positive real number $K$ such that $|f_p(t)|\leq K e$ for all $p\in \enne$ and $t\in [a,b]$.
For every $p\in \mathbb{N}$ and $t \in [a,b]$
we have $|f(t)-f_p(t)|\leq \sigma_p\ e$, and so the  sequence $(f_p)_p$ converges uniformly, and then 
also in $\mu$-measure to $f$. Since $[a,b]$ has finite measure, and 
\[ \int_B |f_p (t)| dt \leq \lambda(B) K\ e,\]
for every Borel set $B\subset [a,b]$ and every integer $p$, the sequence $(f_p)_p$  is equi-absolutely continuous. Moreover, since  
$$\int_a^b|f_{p+r}(t)-f_p(t)|\ dt\leq (b-a)  \sigma_p\ e$$ for all $p$, $r \in \mathbb{N}$, 
we can obtain a mapping $l:\Sigma \to {\bf X}$ as  in Definition \ref{integrabilitainfinita}. Indeed, if 
for all $A \in \Sigma$ we set
$$l_1(A):=\liminf_p \int_A f_p(t)dt,\quad  l_2(A):=\limsup_p \int_A f_p(t)dt,$$ 
then we easily get
$$l_2(A)-\int_A f_p(t)dt \leq \sup_r\int_A(f_{p+r} (t)-f_p(t))dt\leq (b-a)\sigma_p \, e$$
and, for all 
$A \in \Sigma$ and  $p \in \mathbb{N}$
\begin{eqnarray*}
\int_A f_p(t)dt-l_1(A)&\leq& \int_A 
f_p(t)dt-\inf_r\int_Af_{p+r}(t)dt\leq 
\sup_r\int_A(f_p(t)-f_{p+r}(t))dt\leq (b-a)\sigma_p \,e.
\end{eqnarray*}
This clearly implies that $l_1(A)=l_2(A):=l(A)$ for 
all $A \in \Sigma$, and moreover
$$\bigvee_{A\in \Sigma}| l(A)-\int_Af_p(t)dt|
\leq (b-a)\sigma_pe$$
for all $p$, which proves formula (\ref{zvezda})
and hence integrability of $f$.\\
Finally, the proof of differentiability of $\displaystyle{F(x)=\int_a^x \, f(t) \, dt}$, $x \in [a,b]$ is straightforward. $\Box$\\

We are now in position to introduce the moment operators.\\

For every  $f \in \widehat{L} (\lambda)$ (see (\ref{integrale-cappuccio})) let us consider
the operators $T_n f: \mathbb{R}^+ \to {\bf X}$ defined by
\begin{eqnarray}\label{tn}
T_n f (s)=\int_0^{+\infty} M_n \Bigl( \frac{t}{s} \Bigr) f(t) \frac{dt}{t},
\quad n \in \mathbb{N}, \, s \in \mathbb{R}^+.
\end{eqnarray}

In what follows, we shall prove that the sequence $(T_n f - f)_n$ satisfies the
hypotheses of Theorem \ref{trediecibis}, and therefore it is modularly
$(o_{\mathcal{F}})$-convergent to $0$.
The proof of this result and of the next Theorem \ref{unico} will be given in Subsection \ref{prova}.
\begin{teorema}\label{unico}

Let $f$ be a function with compact support $C \subset  [a,b]$ \ with $a > 0$
and belonging to  $\widehat{L}(\lambda)$. Then we have the following:
\begin{description}
\item[\rm \textbf{(\ref{unico}.1)}] if $f$ is bounded, and \, $a> 1$, then the functions $T_n f 
$ are Lipschitz and integrable for every $n \geq 2$; 
 moreover, if $f$ is uniformly continuous, we have 
\begin{eqnarray}{\label{consfubini}}
\int_0^{+\infty}T_n f(s)ds=\frac{n}{n-1}\int_a^bf(s)ds;\end{eqnarray}
\item[\rm \textbf{(\ref{unico}.2)}] if $f$ is uniformly continuous, then $(T_n f)_n$ converges uniformly to $f$;\\
\item[\rm \textbf{(\ref{unico}.3)}] if $f$ is uniformly continuous  and $a >1$, then
\begin{description}
\item[\rm \textbf{(\ref{unico}.3.1)}] \quad
\( \displaystyle{(o)\lim_n \int_0^{+\infty} |T_n f (s)  - f(s)| ds = 0};\)
\item[\rm \textbf{(\ref{unico}.3.2)}] \quad for 
every $n \in \mathbb{N}$, $T_n f $ is uniformly
differentiable and its uniform derivative satisfies
the following relation:
\begin{eqnarray*}
\frac{s}{n}(T_nf)'(s)+T_nf(s)=s(T_1 f)'(s)+T_1 f(s)=f(s)
\end{eqnarray*}
\end{description}
\end{description}
\end{teorema}

\begin{remark}\rm
Replacing $n$ with a real positive constant $\nu$, we can see that, also in our general setting, the linear differential equation
\begin{eqnarray}{\label{riferimento}}
s\varphi'(s)+\nu\varphi(s)=\nu f(s)
\end{eqnarray}
(with $f$ uniformly continuous) can be solved by setting
$$\varphi(s)=\frac{1}{s^{\nu}}\left(c+\nu\int_0^sf(t)t^{\nu-1}dt\right),$$
where $c$ is an arbitrary constant value.
\end{remark}

\begin{remark}\rm
A first application can be formulated in terms of weak convergence of the derivative mappings $(T_nf)'$.
Let $f:[0,+\infty[\to \erre$ be any uniformly continuous map, vanishing outside the interval $[a,T]$, 
with $a >1$, and fix any $C^1$ mapping $w:[0,+\infty[\to \erre$, with compact support.
For each uniformly continuous mapping $g:[0,+\infty[\to \erre$, define
$$\rho(g):=\int_0^{+\infty}|g(s)| \, |w'(s)| \,ds.$$
It is easy to see that $\rho$ is a monotone and 
real-valued modular, and that
$\lim_n\rho(T_nf-f)=0$ thanks to (\ref{unico}.3.1).
 From this we  deduce that
$$\lim_n\int_0^{+\infty}T_nf(s)w'(s)ds =\int_0^{+\infty}f(s)w'(s)ds,$$
and also, integrating by parts:
$$\lim_n\int_0^{+\infty}w(s)(T_nf)'(s)ds=\int_0^{+\infty}w(s)\, df(s)$$
(where the last integral is intended in the Riemann-Stieltjes sense).
Thus we can conclude that, though the mappings $(T_nf)'$ do not converge in general, we always have
$$\lim_n\int_0^{+\infty}w(s)d(T_nf(s))=\int_0^{+\infty}w(s)df(s),$$
for all $C^1$ mappings $w$, with compact support.\\
\end{remark}

\subsection{Brownian motion and stochastic processes}\label{bm}

In order to obtain a more concrete application in Stochastic Integration, we shall assume that ${\bf B}:=(B_t)_{0\leq
t<T}$ (with $T<+\infty$) is the standard Brownian Motion defined on a probability space $(\Omega,\Sigma,P)$.
Of course, we can  consider  ${\bf B}$
as a mapping from $[0,T]$ into $L_2(\Omega)$. 
Since $L_2$ has not a strong unit in general, in order to apply the previous theory we shall suitably restrict it: indeed, thanks to the  well-known Maximum Principle (see e.g. \cite{brei}), we see that there is a suitable positive element $Z\in L_2$ such that $|B(t)|\leq Z$ for all $t\in [0,T]$. 
Moreover, we remark that, (see e.g. \cite{nualart} and \cite{garsia}), there exists a positive random variable $W$ in $L_2$ such that 
$$|B(t+h)-B(t)|\leq |h|^{1/4} W$$
holds, as soon as $t,t+h\in [0,T]$. Thus, taking  ${\bf X}$ as the (complete) subspace of $L_2$ 
generated by all elements dominated by some real  multiple of $W+Z$, we see that ${\bf X}$ has a 
strong unit (i.e. $W+Z$), that $\bf B$ is  ${\bf X}$-valued and  is also a uniformly continuous mapping from 
$[0,T]$ to ${\bf X}$, in the sense of our definition.\\

In order to establish the next results, we introduce a definition.
\begin{definizione}{\label{processiY}} \rm
Let $Y:[0,T]\to {\bf X}$ be any stochastic process, predictable with respect to ${\bf B}$, and with continuous trajectories. If we assume that $\displaystyle{\sup_{t\in [0,T]}E(Y(t)^2)=K<+\infty}$, then it is
well-known that $Y$ is integrable in the It\^{o}'s
sense, with respect to ${\bf B}$ (see e.g. \cite{vara}). Any process $Y$ of this kind will be called a {\em regular} process.
The It\^{o} integral of $Y$ will be denoted as usual with $\displaystyle{(I)\int_0^T Y(t) \,dB(t)}$.
\end{definizione}
\begin{proposizione}{\label{integralepiccolo}}
If the process $Y$ is regular, and $\sup_{t\in [0,T]}E(Y(t)^2)=K$, then
$$\left\|(I)\int_0^TY(t)dB(t)\right\|_2^2\leq KT.$$
\end{proposizione}
{\bf Proof}: Let $D$ be any decomposition of $[0,T]$ obtained by means of the points $t_0=0<t_1<...<t_n=T$, and consider the It\^{o} sum:
$$S(Y,D)=\sum_{j=0}^{n-1}Y(t_i)(B(t_{i+1})-B(t_i)).$$
We easily see that $E(S(Y,D))=0$, since the 
increment $B(t_{i+1})-B(t_i)$ is independent of $Y(t_i)$ for all $i$. Similarly, we deduce that
$cov[Y(t_i)(B(t_{i+1})-B(t_i)),Y(t_j)(B(t_{j+1})-B(t_j)]=0$ as soon as $i\neq j$. So,
\begin{eqnarray*}
E(S(Y,D)^2) &=&V(S(Y,D))=\sum_{i=0}^{n-1}V(Y(t_i)(B(t_{i+1})-B(t_i)))=  \\ &=& \sum_{i=0}^{n-1}E([Y(t_i)(B(t_{i+1})-B(t_i))]^2)=
 \sum_{i=0}^{n-1}
E(Y(t_i)^2)E[(B(t_{i+1})-B(t_i))^2]\leq
\\ &\leq&
 K\sum_{i=0}^{n-1}(t_{i+1}-t_i)=KT.
\end{eqnarray*}
Since the sums $S(Y,D)$ are norm-convergent to the It\^{o} integral $(I)\int_0^TY(t)dB(t)$, the assertion is obvious. $\Box$\\

\medskip

The following corollary is an immediate consequence.
\begin{corollario}{\label{primaconvergenza}}
Let $(Y_n)_n$ be any sequence of regular stochastic processes. Assume that the processes $Y_n$ converge to a regular process $Y$ uniformly in $L_2$, i.e. for every real $\varepsilon>0$ an integer $N$ exists, such that
$$\sup_{t\in [0,T]}\|Y_n(t)-Y(t)\|_2\leq \varepsilon$$
for all $n\geq N$.
Then the It\^{o} integrals $\displaystyle{(I)\int_0^T Y_n(t)dB(t)}$
converge in $L_2$ to the It\^{o} integral \mbox{$\displaystyle{(I)\int_0^TY(t)dB(t)}$.}
\end{corollario}

We shall apply these results to the process $f:[0,+\infty[\to {\bf X}$ defined as follows:
$$f(t):=\left\{\begin{array}{ll}
         0 & t\notin [a,T]\\
         (t-T)(B(t)-B(a)),& t\in [a,T]
\end{array}\right.$$
where $a$ and $T$ are fixed positive numbers, $1<a<T$.
(So $f$ is a process similar to the well-known {\em Brownian Bridge}).

\begin{teorema}{\label{lastfinale}}
Let $f(t)=(t-T)(B(t)-B(a))$ be the process defined above. Then $f$ is clearly ${\bf X}$-valued and uniformly continuous, and we have
$$\lim_n \int_0^T T_nf(s)dB(s)=(I)\int_0^T f(s) dB(s).$$
(We remark that the integral on the left-hand side 
is in the Riemann-Stieltjes sense, since the mappings $T_nf$ are more regular than $f$).
\end{teorema}
{\bf Proof}: It is clear that $f$ is predictable; also  $T_nf$ is for all $n$, since
$T_nf(s)=n\int_0^sf(st)dt$, so the values of $T_nf(s)$ depend only on the values of $f(\tau)$ for $\tau\leq s$; then it is sufficient to apply \ref{integralepiccolo} to the sequence $Y_n:=T_nf-f$, which converges uniformly to 0. $\quad \Box$\\

We note that, in the previous results,  
$(o)$-convergence is not mentioned, since we
don't know whether the integrals $\displaystyle{\int_0^T T_nf(s)dB(s)}$
are dominated by an element in $L_2$. However, this problem  can be solved by replacing the sequence
$\displaystyle{n\mapsto \int_0^T T_nf(s)\,dB(s)}$ with the following one:
$$Y_n:=\left\{\begin{array}{ll}
\left(  \int_0^T T_n f(s)dB(s) \right)\wedge M,& {\rm if }\ \ \int_0^T T_nf(s)dB(s)>0\\
& \\
 \left( \int_0^T T_n f(s)dB(s) \right)\vee (-M),& {\rm if } \ \int_0^T T_nf(s)dB(s)\leq 0,\end{array}\right.$$
 where $M \in {\bf X}$ is any  positive mapping larger than
$\displaystyle{\left|(I)\int_0^T f(s)\,dB(s)\right|}$.

Indeed, $(Y_n)_n$ is obviously dominated, and we always have
$$\left|Y_n-(I)\int_0^T f(s)dB(s)\right|\leq \left|\int_0^T T_nf(s)dB(s)) - (I)\int_0^T f(s)dB(s))\right|,$$
which proves $L_2$-convergence of $(Y_n)_n$ to  the It\^{o} integral of $f$ with respect to ${\bf B}$.  
$\quad \Box$
\\
The following pictures show a comparison among some trajectories of the mapping $f$ and the corresponding ones of the mappings $T_n f$.

\begin{figure}
        \centering     
 \begin{subfigure}[hb]{0.48\textwidth}
                \centering
                \includegraphics[width=\textwidth]{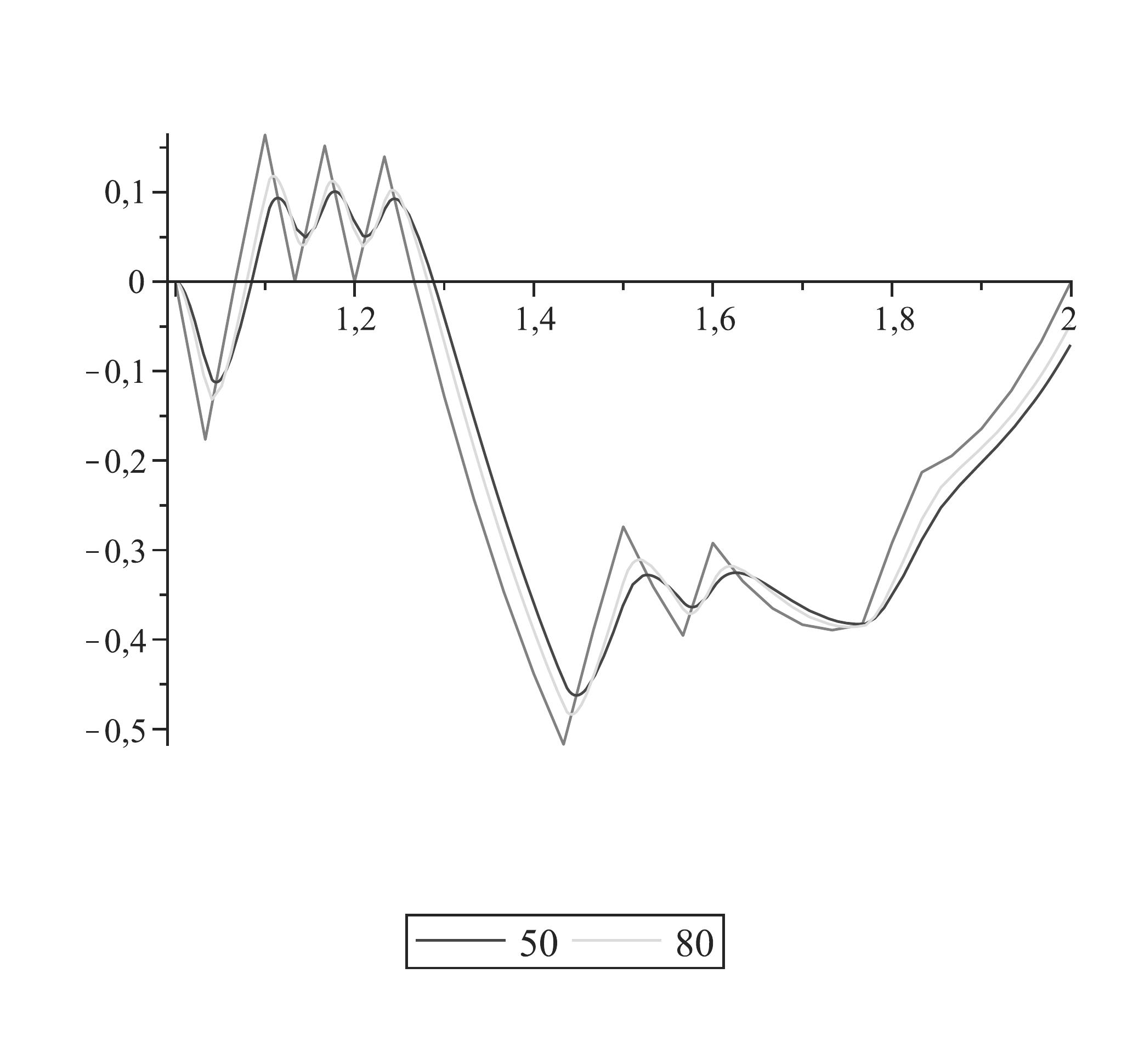}
                \caption{$f$ (gray), $T_{50} f$ (black), $T_{80} f$ (light gray)}
                \label{fig:1}
        \end{subfigure}%
        ~ 
        \begin{subfigure}[hb]{0.48\textwidth}
                \centering
                \includegraphics[width=\textwidth]{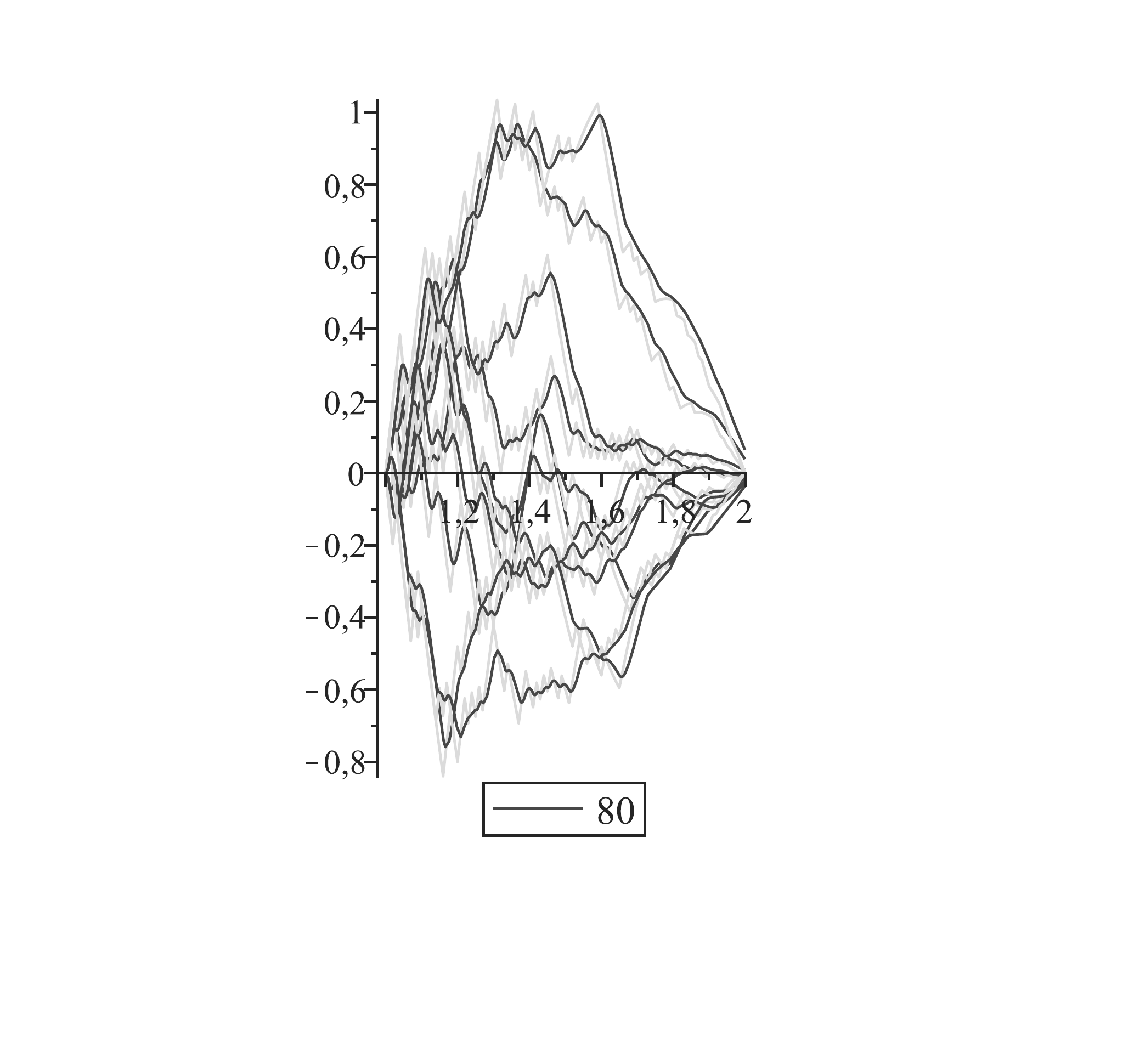}
                \caption{$f$ (ten trajectories  in gray) and  $T_{80} f$ (black)}
                \label{fig:2}
        \end{subfigure}    
        \caption{Pictures of a Brownian Bridge}\label{fig:3}
\end{figure}

\section{Proof of Theorem  \ref{unico}}\label{prova}
Let  $M_f$ denote a fixed real positive number  such that  $\bigvee_{t \in [a,b]} |f(t)|\leq M_f e$  and set $M=n\,b^n M_f$.\\
{\bf Proof of (\ref{unico}.1)}
We first show that $T_n f$ is Lipschitz  and then we will prove integrability of $T_n f$ on  ${\mathbb{R}}^+$.
For each $s \in \mathbb{R}^+$ we can write:
\begin{eqnarray}\label{leggeTnf}
T_n f (s) &=&
n \int_0^{+\infty} \Bigl( \frac{t}{s} \Bigr)^n 1_{(0,1)}
 \Bigl( \frac{t}{s} \Bigr) \frac{f(t)}{t} \, dt= \frac{n}{s^n}
 \int_0^s t^{n-1} f(t) \, dt .
\end{eqnarray}
We now estimate the quantity $|T_n f(s_1) - T_n f(s_2)|$.
Let $s_1$, $s_2 \in {\mathbb R}^+$, with $s_1 <  s_2$. If $s_2 \leq a$, then $T_n f(s_1)-T_n f (s_2)=0$.\\
If $s_2 > a$, then without loss of generality we can suppose $s_1 \geq a$. We get:
\begin{eqnarray*}
T_n f(s_1) - T_n f(s_2)
&=& \frac{n}{s_1^n} \int_0^{s_1} t^{n-1} f(t) \, dt -
\frac{n}{s_2^n} \int_0^{s_1} t^{n-1} f(t) \, dt - \frac{n}{s_2^n} \int_{s_1}^{s_2} t^{n-1} f(t) \, dt=
\\ &=& n \Bigl[   \Bigl( \frac{1}{s_1^n} - \frac{1}
{s_2^n}  \Bigr) \int_0^{s_1} t^{n-1} f(t) \, dt \,
 - \frac{1}{s_2^n} \int_{s_1}^{s_2} t^{n-1} f(t) 
\, dt \Bigr] .
\end{eqnarray*}
Then
\begin{eqnarray*}
|T_n f(s_1)-T_n f(s_2)|\leq M\ \left(\frac{1}
{s_1^n}-\frac{1}{s_2^n}\right) e+M(s_2-s_1)e\leq
M(s_2-s_1)(1+\frac{n}{\tau^{n+1}})e
\end{eqnarray*}
for a suitable $\tau$ with $1 < a \leq s_1 < \tau < s_2$, thanks to the Lagrange Theorem
applied to  the real function $s\mapsto s^{-n}$. So, we see that
$|T_n f(s_1)-T_n f(s_2)|\leq M(1+n)(s_2-s_1)e,$ thus proving that $T_nf$ is Lipschitz.\\
We now prove integrability. Let $n \geq 2$ be fixed. First of all observe
that for every $s>b$ we have
$\displaystyle{T_n \,f(s)= \frac{n}{s^n}K_n}$,
where $K_n\in {\bf X}$ is the integral of $t\mapsto t^{n-1}f(t)$ in $[0,b]$. Since $n\geq 2$, the real-valued mapping $s\mapsto\frac{n}{s^n}$ is clearly integrable in $[b,+\infty[$, and therefore it is easy to see that  $T_n \,f$ is, in the same halfline, thanks to Remark \ref{perdopo}.
Moreover, since  $T_n \,f$ is uniformly continuous, it is also integrable in $[0,b]$, in view of (\ref{uc-prodotto}.4). Therefore we get integrability
 on the whole of $[0,+\infty[$ (see Remark \ref{additivitaestesa}).
Let us now turn to prove (\ref{consfubini}), in the case that $f$ is uniformly continuous.
In order to compute the integral of $T_n \,f$, $n\geq 2$,  we first point out that 
\begin{eqnarray}{\label{kn}}
\int_b^{+\infty} T_n \,f(s)\ ds=K_n \int_b^{+\infty}\frac{n}{s^n}ds=\frac{n}{(n-1)b^{n-1}}K_n, \end{eqnarray}
where $K_n$ is the constant above. 
Now we evaluate the integral in $[a,b]$. To this aim, 
we can observe that, by uniform continuity of $f$ and of the mapping $t\mapsto nt^{n-1}f(t)$ in $[a,b]$, there exist an $(o)$-sequence $(\sigma_p)_p$ in $\erre^+$ and a strictly increasing sequence $(H_p)_p$ in $\enne$, such that 
$$|f(t)-f(u)|+|nt^{n-1}f(t)-nu^{n-1}f(u)|<\sigma_p e$$ 
for all $p\in \enne$ and  $t,u\in [a,b]$ satisfying $|t-u|\leq (b-a)/H_p$. \\
Now, fix $p\in \enne$ and let
$a:=s_0<s_1<s_2<\ldots<s_{H_p}:=b$ be the division  of $[a,b]$ obtained with precisely $H_p$ 
subintervals of the same length, i.e. $(b-a)/H_p$. 
By additivity of the integral we have
$$\int_a^{b} T_n \,f(s) ds=\sum_{j=0}^{H_p-1}\int_{s_j}^{s_{j+1}}T_n \,f(s) ds=\sum_{j=0}^{H_p-1}\int_{s_j}^{s_{j+1}}\frac{1}{s^n}\left[\int_a^s nt^{n-1}f(t)dt\right] ds.$$
For each fixed $s>a$ let us denote by $j(s)$ the maximum index such that $s_{j(s)}< s$, and set also $j(s_0)=j(a)=0$. We have then:
\begin{eqnarray*}
&& \left|  \int_a^snt^{n-1}f(t)dt-
\sum_{j=0}^{j(s)}ns_j^{n-1}f(s_j)\frac{b-a}{H_p}\right| 
\leq  \int_{s_{j(s)}}^snt^{n-1}|f(t)|dt+
n s_{j(s)}^{n-1}|f(s_{j(s)})|\frac{b-a}{H_p} +
\\ &+&
\sum_{j=0}^{j(s)-1}\int_{s_j}^{s_{j+1}}|nt^{n-1}f(t)-ns_{j}^{n-1}f(s_j)|dt
\leq 2 nb^{n-1}M_fe\frac{b-a}{H_p}+\sum_{j=0}^{j(s)-1}\sigma_p (s_{j+1}-s_j)e=
\\ &=&
\left(\frac{2 nb^nM_f}{H_p}+s_{j(s)}\sigma_p \right) e \leq \left(
\frac{2 nb^nM_f}{H_p}+b\sigma_p\right) e.
\end{eqnarray*}
Thus we see that
\begin{eqnarray*}
&|&\int_a^{b} T_n \,f(s) ds-\sum_{j=0}^{H_p-1}\int_{s_j}^{s_{j+1}}\frac{1}{s^n}\sum_{j=0}^{j(s)}ns_j^{n-1}f(s_j)\frac{b-a}{H_p} ds|= \\
&=& |\sum_{j=0}^{H_p-1}\int_{s_j}^{s_{j+1}}\frac{1}{s^n}\left[\int_a^s nt^{n-1}f(t)dt-\sum_{j=0}^{j(s)}ns_j^{n-1}f(s_j)\frac{b-a}{H_p}\right] ds|\leq 
\\ &\leq& \sum_{j=0}^{H_p-1}\int_{s_j}^{s_{j+1}}\frac{1}{s^n}  (
\frac{ 2 nb^nM_f}{H_p}+b\sigma_p)e ds=
 (\frac{ 2 nb^nM_f}{H_p}+b\sigma_p)e\int_a^{b}\frac{1}{s^{n}}ds= \\
&=& \frac{b^{n-1}-a^{n-1}}{n-1} (
\frac{2 nb^nM_f}{H_p}+ b\sigma_p)e.
\end{eqnarray*}
 Letting $\displaystyle{\lambda_p:=\frac{b^{n-1}-a^{n-1}}{n-1} (
\frac{2 nb^nM_f}{H_p}+ b\sigma_p)}$,
 we can see that $(\lambda_p)_p$ is an $(o)$-sequence in $\erre$, and
 summarize the previous result in the following way:
 \begin{eqnarray}{\label{appro1}}
\left|\int_a^{b} T_n \,f(s) 
ds-\sum_{j=0}^{H_p-1}\int_{s_j}^{s_{j+1}}\frac{1}
{s^n}\sum_{j=0}^{j(s)}ns_j^{n-1}f(s_j)\frac{b-
a}{H_p} ds\right|\leq \lambda_p e.\end{eqnarray}
Now, by changing order in the summations, we get
\begin{eqnarray*}
&& \sum_{j=0}^{H_p-1} \int_{s_j}^{s_{j+1}}\frac{1}{s^n}\sum_{j=0}^{j(s)}ns_j^{n-1}f(s_j) 
\frac{b-a}{H_p}ds= 
 \frac{n}{n-1}\frac{b-a}{H_p} \left(\int_{s_0}^{s_1}(n-1)s^{-n}s_0^{n-1}f(s_0) ds+ \right.
\\ &+& \left.
\int_{s_1}^{s_2}[(n-1)s^{-n}s_0^{n-1}f(s_0)+(n-1)s^{-n}s_1^{n-1}f(s_1)] ds+...\right)=
\\ &=& \frac{n}{n-1}\frac{b-a}{H_p}
\left( \int_{s_0}^{b}(n-1)s^{-n}s_0^{n-1}f(s_0)ds+\int_{s_1}^{b}(n-1)s^{-n}s_1^{n-1}f(s_1)+ \ldots + \right.
\\ 
&+& \left. \int_{s_{H_p-1}}^{b}(n-1)s^{-n}s_{{H_p-1}}^{n-1}f(s_{H_p-1}) ds \right)=
\frac{n}{n-1}\frac{b-a}{H_p} \left (\left[ f(s_0)-\frac{s_0^{n-1}f(s_0)}{b^{n-1}}\right]+ \right.
\\ 
\end{eqnarray*}
\begin{eqnarray*}
&+&  \left.
\left[f(s_1)-\frac{s_1^{n-1}f(s_1)}{b^{n-1}}\right]+ \ldots + 
\left[f(s_{H_p-1})-\frac{s_{H_p-1}^{n-1}f(s_{H_p-1})}{b^{n-1}}\right] \right)=\\
&=&\frac{n}{n-1}\sum_{j=0}^{H_p-1}f(s_j)(s_{j+1}-s_j)-\frac{n}{n-1}\sum_{j=0}^{H_p-1}\frac{s_j^{n-1}f(s_j)}{b^{n-1}}(s_{j+1}-s_j).
\end{eqnarray*}

Now, since
\begin{eqnarray*}
\frac{n}{n-1} \left(\left|\sum_{j=0}^{H_p-1}f(s_j)(s_{j+1}-s_j)-\int_a^b f(s) ds\right| \right) &\leq&  \frac{n}{n-1}\sum_{j=0}^{H_p-1}\sigma_p (s_{j+1}-s_j)e =
\frac{n}{n-1}(b-a)\sigma_p e, 
\end{eqnarray*}
and
\begin{eqnarray*}
 && \left|\ \frac{n}{n-1}\sum_{j=0}^{H_p-1}\frac{s_j^{n-1}f(s_j)}{b^{n-1}}(s_{j+1}-s_j)-\frac{n}{(n-1)b^{n-1}}\int_a^bs^{n-1}f(s)ds\right|\leq \\
&&\leq \frac{n}{(n-1)b^{n-1}}\sum_{j=0}^{H_p-1}\sigma_p(s_{j+1}-s_j)e\leq \frac{n}{(n-1)b^{n-2}}\sigma_p e,
\end{eqnarray*}
we obtain
\begin{eqnarray}{\label{appro2}}
\left|\sum_{j=0}^{H_p-1}\int_{s_j}^{s_{j+1}}\frac{1}{s^n}\sum_{j=0}^{j(s)}ns_j^{n-1}f(s_j)\frac{b-a}{H_p} ds- \frac{n}{n-1}\left( \int_a^bf(s)ds- \frac{1}{b^{n-1}}\int_a^bs^{n-1}f(s)ds \right)\right|\leq
h\sigma_p e, \,\,\,
\end{eqnarray}
where $h=2(b-a)+2$.
In conclusion, from (\ref{kn}), (\ref{appro1}) and (\ref{appro2}) we deduce that
\begin{eqnarray*}
&& \left|\int_a^{+\infty}T_nf(s)ds-\frac{n}{n-1}\int_a^bf(s)ds\right|=\\
&=&\left|\int_a^{b}T_nf(s)ds-\frac{n}{n-1}\int_a^bf(s)ds+\frac{n}{(n-1)b^{n-1}}K_n
\right|\leq \lambda_pe+h\sigma_pe
\end{eqnarray*}
holds, for all $p\in \enne$. Letting $p$ tend to $+\infty$, we obtain finally
$$\left|\int_a^{+\infty}T_nf(s)ds-\frac{n}{n-1}\int_a^bf(s)ds\right|=0.\quad \Box$$

\vskip.2cm \noindent
{\bf Proof of (\ref{unico}.2)}
 Let $f$ be as in the hypothesis.
For every $n \in \mathbb{N}$ and $s \in \mathbb{R}^+$, we have
$$|T_n \,f(s)-f(s)|=|\frac{n}{s^n}\int_0^s t^{n-1}f(t)dt-\frac{n}{s^n}\int_0^s t^{n-1}f(s)dt|\leq
\frac{n}{s^n}\int_0^s t^{n-1}|f(t)-f(s)| dt.$$
Since $f$ is uniformly continuous, there exist two $(o)$-sequences $(\sigma_p)_p$ and $(\gamma_p)_p$ in $]0,1[$ such that
$|f(u)-f(v)|\leq  \sigma_p e$ as soon as $|u-v|\leq \gamma_p$ holds, for all positive integers $p$. Now, for every $p\in \enne$ we get
\begin{eqnarray*}
|T_n \,f(s)-f(s)| &\leq& \frac{n}{s^n}\int_0^{s-\gamma_p} t^{n-1}|f(t)-f(s)| dt+\frac{n}{s^n}\int_{s-\gamma_p}^s t^{n-1}|f(t)-f(s)| dt \leq\\
&\leq& \frac{2M_f}{s^n}(s-\gamma_p)^n e+\frac{n}{s^n}\sigma_p\ e\int_0^s\ t^{n-1}dt=2M_f(1-\frac{\gamma_p}{s})^n e+\sigma_p e.
\end{eqnarray*}
Now, for all $s\geq 2b$ we see that
$$|T_n \,f(s)-f(s)|=|T_n \,f(s)|=|\frac{n} {s^n}\int_0^bt^{n-1}f(t)dt|\leq M_f(\frac{b}
{s})^n e\leq \frac{M_f}{2^n} e,$$
while for $s\leq 2b$ we get
$$|T_n \,f(s)-f(s)|\leq 2M_f (1-\frac{\gamma_p}{2b})^n e+\sigma_p e;$$
in any case, we find:
$$\bigvee_{s>0}|T_n \,f(s)-f(s)|\leq 2M_f(1-\frac{\gamma_p}{2b})^n e+\sigma_p\ e +\frac{M_f}{2^n}e,$$
hence
$$\limsup_n\ \bigvee_{s>0}|T_n \,f(s)-f(s)|\leq \sigma_p e$$
holds, for all $p$: this obviously implies that $\limsup_n\ \bigvee_{s>0}|T_n \,f(s)-f(s)|=0$.$\quad\Box$\\

\noindent \bf{Proof of (\ref{unico}.3)}
\rm
\begin{description}
\item[\bf{ (\ref{unico}.3.1)}]
By  (\ref{unico}.1) $T_n f \in L (\lambda)$  for every $n \in \mathbb{N}$, by (\ref{unico}.2)  $(T_n f - f)_n$ $\mu$-converges to zero
and so it is enough to prove equi-absolute continuity of $(T_n f - f)_n$
with respect to the modular $\rho(\cdot)=\int_{ \mathbb{R}^+} |\cdot|\, ds$,
in order to apply Theorem \ref{trediecibis}.
From uniform convergence of $T_n f$ to $f$ we see that, setting
$ \zeta_n := \sup_{p\geq n}\bigvee_{s \in \mathbb{R}^+} | T_p f(s) - f(s)|,$
$(\zeta_n)_n$ is an $(o)$-sequence. Then
\begin{eqnarray*}
\rho([T_n f - f]\, 1_B) &=&  \int_B \Bigl|T_n  f (s) - f(s)\Bigr| \, ds  \leq
\zeta_n  \lambda(B)
\end{eqnarray*}
for every $n$ and whenever  $B \in \Sigma$.
This clearly implies $ac_{\rho}(1)$.
We now turn to ($ac_{\rho}(2)$). For each $m \in \mathbb{N}$, let $B_m=[0, 2b]$.
As we already pointed out, when $s>2b$ it holds: 
$|T_n \,f(s)|\leq M_f (b/s)^n e.$ 
So, for $m$, $n \in \mathbb{N}$ we have
\begin{eqnarray*}
\rho([T_n f - f] \,1_{\mathbb{R}^+ \setminus B_m})  \leq
M_f b^n e\int_{2b}^{+\infty} \frac{1}{s^n} ds  \leq M_ f \ b^n \ \frac{1}{(2b)^{n-1}}e= \frac{M_f b}{2^{n-1}}e.
\end{eqnarray*}
This clearly means that \( (r)\lim_n \rho([T_n f  - f] 
\,1_{\mathbb{R}^+ \setminus B_m}) =0 \) 
uniformly with respect to $m$, and hence 
($ac_{\rho}(2)$).
So, all hypotheses of the Vitali theorem are 
satisfied and the assertion
follows from Theorem \ref{trediecibis}.
$\quad\Box$\\
\item[\bf{ (\ref{unico}.3.2)}]
We shall use the following expression for $T_n f$:
\begin{eqnarray}\label{seconda}
T_n f (s) = \frac{n}{s^n}
 \int_0^s t^{n-1} f(t) \, dt .
\end{eqnarray}
We first observe that the integrand  
$(\cdot)^{n-1} \, f(\cdot):[0,s] \to {\bf X}$
satisfies (\ref{uc-prodotto}.1) and hence its
integral function is uniformly differentiable by (\ref{uc-prodotto}.4).
So $T_n f$  is the product of two functions which satisfy  (\ref{uc-prodotto}.3), i.e.  $n \, s^{-n}$ and
$s\mapsto \int_0^s t^{n-1} \, f(t)$.
By Proposition \ref{uc-prodotto}
we obtain:
\begin{eqnarray}{\label{recupero}}
\left( T_n f \right)^{\prime} (s) = -\frac{n^2}{s^{n+1}} \int_0^s t^{n-1} f(t) dt + \frac{n}{s} f(s)
=
\frac{n}
{s}\big(f(s)-T_n(f)(s) \big).
\end{eqnarray}
Formula (\ref{recupero}) shows that $f$ is 
uniquely determined by each of the functions $T_n 
f$'s, in particular
\begin{eqnarray*}
f(s)=s(T_1 f)'(s)+T_1 f(s)
\end{eqnarray*}
holds for all $s$, and therefore all functions $T_n 
f$'s can be easily obtained from one another.
$\quad\Box$\\
\end{description}
\bibliographystyle{elsarticle-num}

\end{document}